\documentclass[12pt]{amsart}
\usepackage{amsmath}
\usepackage{amsfonts}
\usepackage{amsthm}
\usepackage{amssymb}

 \newtheorem{thm}{Theorem}[section]

 \theoremstyle{definition}
 
  \newtheorem{ques}[thm]{Question}
 \theoremstyle{remark}

 \numberwithin{equation}{section}

\newcommand{\D}{\mathbb D}
\newcommand{\C}{\mathbb C}
\newcommand{\hol}{\mathcal Hol}
\newcommand{\ig}{\stackrel{\text{def}}{=}}
\DeclareMathOperator{\og}{O} \DeclareMathOperator{\op}{o}

\DeclareMathOperator{\Real}{Re} 

\newenvironment{Pf}{\noindent{\emph{Proof of}}}{$\Box$}

\newtheorem{other}{\bf Theorem}              

\begin{document}

%
%
%
%
%
%
%
%
%

\title[Products of unbounded Bloch functions]
 {Products of unbounded Bloch functions}

\author[Daniel Girela]{Daniel Girela}

\address{%
An\'{a}lisis Matem\'{a}tico\\
Facultad de Ciencias\\
Universidad de M\'{a}laga\\
29071 M\'{a}laga\\
Spain}

\email{girela@uma.es}

\thanks{This research is supported in part by a grant from \lq\lq El Ministerio de
Econom\'{\i}a y Competitividad\rq\rq , Spain (PGC2018-096166-B-I00)
and by grants from la Junta de Andaluc\'{\i}a (FQM-210 and
UMA18-FEDERJA-002).}



\date{October 25, 2019}

\begin{abstract} We give new constructions of pair of  functions $(f, g)$, analytic
in the unit disc, with $g\in H^\infty $ and $f$ an unbounded Bloch
function, such that the product $g\cdot f$ is not a Bloch function.

\end{abstract}

\maketitle
\section{Introduction and statements of the results}\label{intro}
Let $\D=\{z\in\C:|z|<1\}$ denote the open unit disc in the complex
plane $\mathbb C$. The space of all analytic functions in $\D $ will
be denoted by $\hol (\D )$.
\par\medskip
For $0<p\le\infty $, the classical Hardy space $H\sp p $ is defined
as the set of all $f\in\hol (\D )$ for which
$$\|f\|_{H\sp
p}\ig \sup_{0<r<1}M\sb p(r,f)<\infty ,$$ where, for $0<r<1$ and
$f\in \hol (\D )$,
\begin{align*}M_p(r,f)\,=\,&\left (\frac{1}{2\pi}\int_0^{2\pi}
|f(re\sp{i\theta })|\sp p d\theta\right )\sp{1/p}, \,(0<p<\infty );
\\M_\infty (r,f)\,=\,&\sup _{\theta \in \mathbb R}|f(re\sp{i\theta
})|.
\end{align*} We mention \cite{Du:Hp} as a general reference for
the theory of Hardy spaces.
\par\medskip
A function $f\in\hol (\D )$ is said to be a Bloch function if
$$\Vert f\Vert _{\mathcal B}\ig \vert f(0)\vert +
\sup _ {z\in \D}\,(1-\vert z\vert\sp 2)\vert f\sp\prime (z)\vert
<\infty .
$$
The space of all Bloch functions is denoted by $\mathcal B$, it is a
Banach space with the just defined norm $\Vert \cdot \Vert
_{\mathcal B}$. It is well known that
$$H\sp\infty\subsetneq\mathcal B.$$ A typical example of an
unbounded Bloch function is the function $f$ defined by
$$f(z)\,=\,\log \frac{1}{1-z},\quad z\in\mathbb D.$$
We mention \cite{ACP} as a general reference for the theory of Bloch
functions.
\par\medskip
A function $f$ which is meromorphic in $\D $ is said to be a normal
function in the sense of Lehto and Virtanen \cite{LV} if
$$\sup _{z\in \D }(1-\vert z\vert \sp 2)\frac{\vert f\sp\prime (z)\vert }
{1+\vert f(z)\vert \sp 2}<\infty .$$ For simplicity, we shall let
$\mathcal N $ denote the set of all holomorphic normal functions in
$\D $. It is clear that any Bloch function is a normal function,
that is, we have $\mathcal B\subset \mathcal N$. We refer to
\cite{ACP}, \cite{LV} and \cite{Po:UF} for the theory of normal
 functions. In particular, we remark here that if $f\in \mathcal N$, $\xi \in
\partial \D$ and $f$ has the asymptotic value $L$ at $\xi$, (that
is, there exists a curve $\gamma $ in $\D$ ending at $\xi $ such
that $f(z)\to L$, as $z\to \xi $ along $\gamma $) then $f$ has the
non-tangential limit $L$ at $\xi $.
\par\medskip
Let us recall that if a sequence of points $\{a_n\} $ in  the unit
disc satisfies the \emph{Blaschke condition\/}:
$\sum\sb{n=1}\sp{\infty} (1-|a\sb n|)<\infty $, the corresponding
Blaschke product $B $ is defined as
\[
 B(z)=\prod\sb{n=1}\sp{\infty} \frac{|a\sb{n}|}{a\sb{n}}
 \frac{a\sb{n}-z}{1-\overline{a\sb{n}}z} \,.
\]
Such a product is analytic in $\D $. In fact, it is an inner
function, that is, an $H^\infty $-function with radial limit of
absolute value $1$ at almost every point of $\partial \mathbb D$
(cf. \cite[Chapter~2]{Du:Hp}).
\par If $\{a_n\} $ is a Blaschke sequence and there exists $\delta >0$
such that $$\prod_{m\neq n} \left\vert \frac{a_n-a_m}{1-\overline
a_na_m}\right \vert\ge \delta ,\quad\text{for all $n$,}$$ we say
that the sequence $\{a_n\} $ is \emph{uniformly separated} and that
$B$ is an \emph{interpolating Blaschke product}. Equivalently,
\begin{equation}\label{inter}
B \,\, {\it is \, an \, interpolating\, Blaschke\,  product} \,\,\,
\Leftrightarrow\,\,\, \inf _{n\ge 1}(1-\vert a_n\vert \sp 2)\vert
B'(a_n)\vert >0.\end{equation} We refer to \cite[Chapter~9]{Du:Hp}
and \cite[Chapter~VII]{Gar} for the basic properties of
interpolating Blaschke products. In particular, we recall that an
exponential sequence is uniformly separated and that the converse
holds if all the $a_k$'s are positive.
\par\medskip
Lappan \cite[Theorem~3]{La} proved that if $B$ is an interpolating
Blaschke product and $f$ is a normal analytic function in $\D $, the
product $B\cdot f$ need not be normal. Lappan used this to show that
$\mathcal N$ is not a vector space.
\par
Lappan's result is a consequence of the following easy fact: if $B$
is an interpolating Blaschke product whose sequence of zeros is $\{
a_n\} $ and $G$ is an analytic function in $\mathbb D$ with
$G(a_n)\to\infty$, then $f=B\cdot G$ is not a normal function (and
hence it is not a Bloch function either). This result has been used
by several authors (see \cite{Cam,Ya1,Ya2,Gi:Pri,Gi:Ill,BGM}) to
construct distinct classes of non-normal functions.
\par The author and Su\'{a}rez proved in \cite{GS} a result of this
kind dealing with Blaschke products with zeros in a Stolz angle but
not necessarily interpolating, improving a result of \cite{GGP}.
Namely, Theorem\,\@1 of \cite{GS} is the following.
\begin{other}\label{Th-GGP}
Let $B$ be an infinite Blaschke product whose sequence of zeros $\{
a_n\} $ is contained in a Stolz angle with vertex at $1$ and let $G$
be analytic in $\mathbb D$ with $G(z)\to\infty$, as $z\to 1$. Then
the function $f=B\cdot G$ is not a normal function.\end{other}
\par\medskip It is natural to ask whether it is possible to
prove results similar to those described, substituting \lq\lq
Blaschke products\rq\rq \, by some other classes of $H^\infty
$-functions. Our first result in this paper deals with the atomic
singular inner function.
\begin{thm}\label{Bloch-inner} Let $S$ be the atomic singular inner
function defined by $$S(z)=\exp\left (-\frac{1+z}{1-z}\right ),\quad
z\in \mathbb D,$$ and let $f$ be a Bloch function with
$$\lim_{z\to 1}\vert f(z)\vert =\infty .$$
Then the function $F$ defined by $F(z)=S(z)f(z)$ is not a normal
function (hence, it is not a Bloch function).\end{thm}
\par\medskip
In particular, the function $f$ defined by $f(z)=S(z)\cdot \log
\frac{1}{1-z}$ ($z\in \mathbb D$) is not normal. Since a Bloch
function satisfies $M_\infty (r,f)\,=\,\og \left
(\log\frac{1}{1-r}\right )$, Theorem\,\@\ref{Bloch-inner} follows
from the following result.

\begin{thm}\label{unbound-inner}
Let $S$ be the singular inner function defined by $S(z)=\exp\left
(-\frac{1+z}{1-z}\right )$ ($z\in \mathbb D$) and let $f$ be an
analytic function in $\mathbb D$ satisfying:
\begin{itemize}\item[(i)] $\lim_{z\to 1}\vert f(z)\vert =\infty
$.\item[(ii)] $\vert f(r)\vert \,=\,\op \left
(\exp\frac{1+r}{1-r}\right )$,\,as $r\to 1^-$.
\end{itemize}
Then the function $F$ defined by $F(z)=S(z)f(z)$ is not a normal
function (hence, it is not a Bloch function).\end{thm}
\par\medskip
For $g\in\hol (\D )$, the multiplication operator $M_g$ is defined
by \[ M_g(f)(z)\ig g(z)f(z),\quad f\in \hol (\D),\,\, z\in \D.
\]
Let us recall that if $X$ and $Y$ are two spaces of analytic
function in $\D $ and $g\in\hol (\D )$ then $g$ is said to be a
multiplier from $X$ to $Y$ if $M_g(X)\subset Y$. The space of all
multipliers from $X$ to $Y$ will be denoted by $M(X,Y)$ and $M(X)$
will stand for $M(X,X)$. Brown and Shields \cite{BS} characterized
the space of multipliers of the Bloch space $M\mathcal (B)$ as
follows.
\begin{other}\label{mult-Bloch}A function $g\in\hol (\mathbb D)$ is
a multiplier of the Bloch space if and only if $g\in H^\infty \cap
 \mathcal B_{\log }$, where $\mathcal B_{\log }$ is the Banach space of those
functions $f\in\hol (\D )$ which satisfy
\begin{equation}\label{Blog}\Vert f\Vert _{B_{\log }}\ig
\vert f(0)\vert +\sup_{z\in \mathbb D}(1-\vert z\vert ^2)\left (\log
\frac{2}{1-\vert z\vert ^2}\right )\vert f^\prime (z)\vert <\infty
.\end{equation}
\end{other}
Thus, if $g\in H^\infty \setminus \mathcal B_{\log }$ there exists a
function $f\in \mathcal B\setminus H^\infty $ such that $g\cdot
f\notin \mathcal B$. It is easy to see that the analytic Lipschitz
spaces $\Lambda _\alpha $ ($0<\alpha \le 1$) and the mean Lipschitz
spaces $\Lambda ^p_{\alpha }$ ($1<p<\infty $, $1/p<\alpha \le 1$)
are contained in $M(\mathcal B)$ We refer to
\cite[Chapter\,\@5]{Du:Hp} for the definitions of these spaces, let
us simply recall here that
$$\Lambda ^1_1=\{ f\in \hol (\mathbb D) : f^\prime \in H^1\} .$$
On the other hand, Theorem\,\@1 of \cite{GGH} shows the existence of
a Jordan domain $\Omega $ with rectifiable boundary and $0\in \Omega
$, and such that the conformal mapping $g$ from $\mathbb D$ onto
$\Omega $ with $g(0)=0$ and $g^\prime (0)>0$ does not belong to
$\mathcal B_{\log }$. For this function $g$ we have that $g\in
\Lambda ^1_1$ but $g$ is  not a multiplier of $\mathcal B$. Thus we
have:
$$\Lambda ^1_1\,\not\subset \,M(\mathcal B).$$
\par In view of this and the results involving Blaschke products that we have
mentioned above, it is natural to ask the following question:
\begin{ques}\label{question} Is it true that for any given $f\in
\mathcal B\setminus H^\infty $ there exists a function $g\in \Lambda
^1_1$ such that $g\cdot f\notin \mathcal B$?
\end{ques}
\par\medskip We shall show that the answer to this question is
affirmative. Actually we shall prove a stronger result.
\par We let $B^1$ denote the minimal Besov space which consists of those
functions $f\in \hol (\mathbb D)$ such that $$\int_{\mathbb D}\vert
f^{\prime \prime }(z)\vert \,dA(z)\,<\,\infty .$$ Here $dA$ denotes
the area measure on $\mathbb D$. Alternatively, the space $B^1$ can
be characterized as follows (see \cite{AFP}): \par For $f\in\hol (\D
)$, we have that $f\in B^1$ if and only there exist a sequence of
points $\{ a_k\} _{k=1}^\infty \subset \D $  and a sequence $\{
\lambda_k\} _{k=0}^\infty \in \ell ^1$ such that
\begin{equation}\label{B1def}
f(z)=\lambda _0+\sum_{k=1}^\infty\lambda_k\varphi_{a_k}(z),\quad
z\in\D.\end{equation} Here, for $a\in \mathbb D$, $\varphi
_a:\mathbb D\rightarrow \mathbb D$ denotes the M\"{o}bius
transformation defined by
\begin{equation}\label{phia}\varphi_a(z)=\frac{a-z}{1-\overline az},\quad z\in \mathbb D
.\end{equation} is is well known that $B^1\subset \Lambda^1_1$ (see
\cite{AFP,DGV1}) and then our next result implies that the answer to
question\,\@\ref{question} is affirmative.

\par\medskip
\begin{thm}\label{not-mul-B1} If $f\in \mathcal B\setminus H^\infty $ then there
exists $g\in B^1$ such that $g\cdot f\notin \mathcal B$.
\end{thm}
\par\medskip The proofs of Theorem\,\@\ref{unbound-inner} and
Theorem\,\@\ref{not-mul-B1} will be presented in
section\,\@\ref{proofs}. We close this section noticing that
throughout the paper we shall be using the convention that $C=C(p,
\alpha ,q,\beta , \dots )$ will denote a positive constant which
depends only upon the displayed parameters $p, \alpha , q, \beta
\dots $ (which often will be omitted) but not necessarily the same
at different occurrences. Moreover, for two real-valued functions
$E_1, E_2$ we write $E_1\lesssim E_2$, or $E_1\gtrsim E_2$, if there
exists a positive constant $C$ independent of the arguments such
that $E_1\leq C E_2$, respectively $E_1\ge C E_2$. If we have
$E_1\lesssim E_2$ and $E_1\gtrsim E_2$ simultaneously then we say
that $E_1$ and $E_2$ are equivalent and we write $E_1\asymp E_2$.

\par\medskip
\section{The proofs}\label{proofs}
\par\medskip
\begin{Pf}{\it Theorem\,\@\ref{unbound-inner}.}
For $0<a<1$, set $\Gamma _a=\{ z\in \mathbb D : \vert z-a\vert
={1-a}\} $. If $z\in \Gamma _a$ then $\Real
\frac{1+z}{1-z}\,=\,\frac{a}{1-a}$ and, hence,
\begin{equation*}\label{modSGammaa}\vert S(z)\vert \,=\,\exp\left
(\frac{-a}{1-a}\right ),\quad z\in \Gamma _a.\end{equation*} This,
together with (i), implies that
$$F(z)\to \infty ,\quad \text{as $z\to 1$ along $\Gamma _a$}.$$
Hence $F$ has the asymptotic value $\infty $ at $1$. On the other
hand, (ii) implies that $F$ has the radial limit $0$ at $1$. Then it
follows that $F$ is not normal.
\end{Pf}
\par\medskip

\begin{Pf}{\it Theorem\,\@\ref{not-mul-B1}.} Take $f\in \mathcal B\setminus H^\infty $. Set $$\varphi (r)\,=\,M_\infty (r,g),\quad 0<r<1.$$
Clearly, $\varphi (r)\to\infty $, as $r\to 1$ and it is well known
that $$\phi (r)\,=\,\og \left (\log \frac{1}{1-r}\right ).$$ This
implies that \begin{equation}\label{op-square}(1-r)^2\varphi (r)\to
0,\quad \text{as $r\to 1$}.\end{equation} Choose a sequence of
numbers $\{ r_n\} \subset (0,1)$ satisfying the following
properties:
\begin{itemize}\item[(i)] $\{ r_n\} $ is increasing.
\item[(ii)] $(1-r_n)^2\varphi (r_n)\,=\,\op \left (\left (
\frac{1-r_{n-1}}{n}\right )^2\right ),\quad \text{as $n\to \infty
$.}$
\item[(iii)] $\varphi (r_n)\,\ge\,2\varphi (r_{n-1})$, for all $n$.
\item[(iv)] $\frac{1-r_{n+1}}{1-r_n}\to 0$, as $n\to \infty $.
\end{itemize}
The existence of such a sequence is clear, bearing in mind
(\ref{op-square}) and the the fact that $\varphi (r)\to\infty $, as
$r\to 1$.
\par Set
$$\lambda _k\,=\,\varphi (r_k)^{-1/2},\quad k=1, 2, \dots .$$
For each $k$, take $a_k\in \mathbb  D$ with $\vert a_k\vert
\,=\,r_k$ and $\vert f(a_k)\vert \,=\,\varphi (r_k)$. Using (iii),
it follows that \begin{equation}\label{sumlam}\sum_{k=1}^\infty
\lambda _k\,<\,\infty . \end{equation} Define
\begin{equation}\label{def-f}g(z)\,=\,\sum_{k=1}^\infty \lambda
_k\varphi _(z),\quad z\in \mathbb D.\end{equation} Using
(\ref{sumlam}) we see that the sum in (\ref{def-f}) defines an
analytic function in $\mathbb D$ which belongs to $B^1$. Set
$$F(z)=g(z)f(z),\quad z\in \mathbb D.$$ Since $g\in H^\infty $ and
$f\in\mathcal B$ we see that
\begin{equation}\label{gfprime}\vert g(a_n)f^\prime (a_n)\vert
\,=\,\og \left (\frac{1}{1-\vert a_n\vert }\right ).\end{equation}
On the other hand,
\begin{equation}\label{gprime}\vert g^\prime (a_n)f(a_n)\vert
\,\gtrsim \,I\,-\,II\,-\,III,\end{equation} where
\begin{align*}
I\,=&\,\vert f(a_n)\vert \lambda_n \vert
\varphi_{a_n}(a_n)\vert,\,\,\,II\,\lesssim \,\vert f(a_n)\vert
\sum_{k=1}^{n-1}\lambda _k\frac{1-\vert a_k\vert ^2}{\vert
1-\overline {a_k}\,a_n\vert ^2},
\\ III\,\lesssim &\,\vert
f(a_n)\vert \sum_{k=n+1}^{\infty }\lambda _k \frac{1-\vert a_k\vert
^2}{\vert 1-\overline {a_k}\,a_n\vert ^2}.
\end{align*}
Clearly,
\begin{equation}\label{I}
I=\,\vert f(a_n)\vert \lambda_n \vert \varphi_{a_n}(a_n)\vert
\,\asymp \,\frac{\varphi (r_n)^{1/2}}{1-r_n}.\end{equation} Using
the definitions, the facts that $\varphi $ and the sequence $\{
r_n\} $ are increasing, and (ii), we obtain
\begin{align}\label{II}II\,\lesssim \,&\vert f(a_n)\vert
\sum_{k=1}^{n-1}\lambda _k\frac{1-\vert a_k\vert ^2}{\vert
1-\overline {a_k}\,a_n\vert ^2}\\ \,\lesssim \,& \varphi
(r_n)\sum_{k=1}^{n-1}\varphi (r_k)^{-1/2}\frac{1-\vert a_k\vert
}{[(1-\vert a_k\vert )+(1-\vert a_n)]^2} \nonumber
\\ \,\lesssim \,&\varphi
(r_n)\sum_{k=1}^{n-1}\frac{1}{\varphi (r_k)^{1/2}(1-r_k)} \nonumber
\\
\,\lesssim \,&\frac{n\varphi (r_n)}{1-r_{n-1}} \nonumber
\\
\,=\,& \frac{\varphi (r_n)^{1/2}}{1-r_n}\,\varphi
(r_n)^{1/2}\,\frac{n(1-r_n)}{1-r_{n-1}}\nonumber
\\
\,=\,&\op \left (\frac{\varphi (r_n)^{1/2}}{1-r_n}\right )\nonumber.
\end{align}
Likewise, using the definitions, the facts that $\varphi $ and the
sequence $\{ r_n\} $ are increasing, (iii), and (iv),  we obtain
\begin{align}\label{III}
III\,\lesssim \,&\varphi (r_n)\sum_{k=n+1}^\infty \frac{\varphi
(r_k)^{-1/2}(1-r_k)}{[(1-r_k)+(1-r_n)]^2}
\\ \,\lesssim \,&\varphi
(r_n)\sum_{k=n+1}^\infty \varphi
(r_k)^{-1/2}\frac{1-r_k}{(1-r_n)^2}\nonumber
\\ \,\lesssim \,&\varphi
(r_n)\frac{1-r_{n+1}}{(1-r_n)^2}\sum_{k=n+1}^\infty \varphi
(r_k)^{-1/2}\nonumber
\\ \,\lesssim \,&
\frac{\varphi (r_n)^{1/2}}{1-r_n}\cdot
\frac{1-r_{n+1}}{1-r_n}\nonumber
\\
\,=\,&\op \left (\frac{\varphi (r_n)^{1/2}}{1-r_n}\right ).\nonumber
\end{align}
Using (\ref{I}), (\ref{II}), (\ref{III}), and the fact that
$\lim\varphi (r_n)=\infty $, we deduce that
$$(1-\vert a_n\vert )\vert g^\prime (a_n)f(a_n)\vert \,\rightarrow
\,\infty ,\quad\text{as $n\to\infty$}.$$ This and (\ref{gfprime})
imply that $F$ is not a Bloch function.
\end{Pf}

\medskip
\bibliographystyle{amsplain}

\end{document}